\documentclass{amsart}
\usepackage{amssymb}
\usepackage{amscd}
\usepackage{amsthm}

\newtheorem{thm}{Theorem}[section]
\newtheorem{lem}[thm]{Lemma}
\newtheorem{cor}[thm]{Corollary}
\newtheorem{prop}[thm]{Proposition}
\theoremstyle{definition}

\newtheorem{rem}[thm]{Remark}

\newtheorem{sect}[thm]{}

\def \kbar {{\bar k}}

\def \Romannumeral #1 {\expandafter\uppercase\expandafter
{\romannumeral #1} }

\def \Br {{\rm{Br}}}

\def \Ga {{\Gamma}}
\def \bR {{\bf R}}
\def \Pic {{\rm {Pic\,}}}

\def \Ker {{\rm{Ker\,}}}
\def \Im {{\rm {Im\,}}}

\def \A{{\mathbb A}}

\def \Spec {{\rm{Spec\,}}}
\def \dim {{\rm{dim\,}}}
\def \Hom {{\rm {Hom}}}

\def \Aut{{\rm Aut \,}}

\def\ov{\overline}

\def\lra{\longrightarrow}

\def \Z {{\mathbb Z}}
\def \Q {{\mathbb Q}}
\def \F {{\mathbb F}}
\def \Ext {{\rm Ext}}

\def \non {{\rm{non-}}}
\def\e{{\varepsilon}}

\def\G{{\bf G}}
\def\R{{\mathbb R}}
\def\P{{\mathbb P}}
\def\D{{\mathcal D}}

\def\T{{\cal T}}
\def\lra{\longrightarrow}

\def\H{{\rm H}}

\def\id{{\rm id}}

\def\sH{{\mathcal H}}
\def\T{{\mathcal T}}

\def\NS{{\rm NS\,}}
\def\Tor{{\rm Tor}}

\def\sF{{\mathcal F}}

\def\Ga{\Gamma}

\def\e{\varepsilon}
\def\et{{\rm \acute et}}
\def\tors{{\rm tors}}

\def\fchar{{\rm char}}

      \def\Spec{{\rm Spec}}
           
\newcommand{\bthe}{\begin{theo}}
\newcommand{\ble}{\begin{lem}}
\newcommand{\bpr}{\begin{prop}}
\newcommand{\bco}{\begin{cor}}
\newcommand{\bde}{\begin{defi}}
\newcommand{\ethe}{\end{theo}}
\newcommand{\ele}{\end{lem}}
\newcommand{\epr}{\end{prop}}
\newcommand{\eco}{\end{cor}}
\newcommand{\ede}{\end{defi}}

\def\bsm{\left( \begin{smallmatrix}}
\def\esm{\end{smallmatrix} \right)}

\title{The Brauer group and the Brauer--Manin set of products of varieties}
\author[Alexei\ N.\ Skorobogatov]{Alexei\ N.\ Skorobogatov}
\address{Department of Mathematics, South Kensington Campus,
Imperial College London, SW7 2BZ England, U.K.}
\address{Institute for the Information Transmission Problems, Russian
Academy of Sciences, 19 Bolshoi Karetnyi, Moscow, 127994 Russia}
\email{a.skorobogatov\char`\@imperial.ac.uk}
\author[Yuri\ G.\ Zarhin]{Yuri\ G.\ Zarhin}
\address{Department of Mathematics, Pennsylvania State University,
University Park, PA 16802, USA}
\address{Institute for Mathematical Problems in Biology,
Russian Academy of Sciences, Pushchino, Moscow Region, Russia}
\email{zarhin\char`\@math.psu.edu}

\begin{document}

\baselineskip=15pt
%\begin{abstract}
%\end{abstract}

\maketitle

\begin{abstract}
Let $X$ and $Y$ be smooth and projective varieties over a field $k$ 
finitely generated over $\Q$, and let $\ov X$ and $\ov Y$ be the 
varieties over an algebraic closure of $k$ obtained
from $X$ and $Y$, respectively, by extension of the ground field.
We show that the Galois invariant
subgroup of $\Br(\ov X)\oplus \Br(\ov Y)$ has finite index
in the Galois invariant subgroup of $\Br(\ov X\times\ov Y)$.
This implies that the cokernel of the natural map
$\Br (X)\oplus\Br (Y)\to\Br(X\times Y)$ is finite when $k$ is a number field.
In this case we prove that the Brauer--Manin set of the product of
varieties is the product of their Brauer--Manin sets.
\end{abstract}

\bigskip

Let $k$ be a field with a separable closure $\bar{k}$, $\Gamma=\Aut(\bar{k}/k)$.
For an algebraic variety $X$ over $k$ we write $\ov X$ for
the variety over $\bar k$ obtained from $X$ by extending the ground field.
Let $\Br(X)$ be the cohomological
Brauer--Grothendieck group $\H^2_\et(X,\G_m)$, see \cite{Gr}.
The group $\Br(\ov X)$ is naturally a Galois module.
The image of the natural homomorphism $\Br(X) \to \Br(\ov X)$ lies in
$\Br(\ov X)^{\Gamma}$; the kernel of this homomorphism is denoted
by $\Br_1(X)$, so that
$\Br(X)/\Br_1(X)$ is a subgroup of $\Br(\ov X)^{\Gamma}$.
Recall that $\Br(X)$ and $\Br(\ov X)$ are torsion abelian groups
whenever $X$ is smooth, see \cite[II, Prop. 1.4]{Gr}.

\medskip

\noindent{\bf Theorem A} {\it
Let $k$ be a field finitely generated over $\Q$. Let $X$ and $Y$ be
smooth, projective and geometrically integral varieties over $k$.
Then the cokernel of the natural injective map
$$\Br(\ov X)^{\Gamma}\oplus\Br(\ov Y)^{\Gamma}\to
\Br(\ov X\times \ov Y)^{\Gamma}$$
is finite.}

\medskip

\noindent See Theorem \ref{main} for an analogue in characteristic $p\not=2$.
The proof uses the results of Faltings and the second named author 
on Tate's conjecture for abelian varieties.

Let $k$ be a field finitely generated over its prime subfield. 
We proved in our previous paper \cite{SZJAG}
that $\Br(\ov X)^\Gamma$ is finite when $X$ is an abelian
variety and $\fchar(k)\ne 2$, or $X$ is a K3 surface and $\fchar(k)=0$.
As a corollary we obtain that if $Z$ is a smooth and projective 
variety over $k$ such that $\ov Z$ is birationally equivalent 
to a product of curves, abelian
varieties and K3 surfaces, then the groups $\Br(\ov Z)^{\Gamma}$
and $\Br(Z)/\Br_1(Z)$ are finite.

The following result easily follows from Theorem A, see Section \ref{3}.

\medskip

\noindent{\bf Theorem B} {\it
Let $k$ be a field finitely generated over $\Q$. Let $X$ and $Y$ be
smooth, projective and geometrically integral varieties over $k$.
Assume that $(X\times Y)(k)\not=\emptyset$ or $\H^3(k,\bar k^*)=0$.
Then the cokernel of the natural map
$$\Br(X)\oplus\Br(Y)\to\Br(X\times Y)$$
is finite.}

Now let $k$ be a number field. In this case $\H^3(k,\bar k^*)=0$, see 
\cite[Cor. I.4.21]{ADT}, so by Theorem B the Brauer group
$\Br(X\times Y)$ is generated, modulo the image of 
$\Br(X)\oplus \Br(Y)$, by finitely many elements. The following result shows
that these elements do not give any new Brauer--Manin
conditions on the adelic points of $X\times Y$
besides those already given by the elements of $\Br(X)\oplus\Br(Y)$.
For the definition of the Brauer--Manin set $X(\A_k)^\Br$ we refer
to \cite[Section 5.2]{S}.

\medskip

\noindent{\bf Theorem C} {\it
Let $X$ and $Y$ be smooth, projective, geometrically integral
varieties over a number field $k$. Then we have
$$(X\times Y)(\A_k)^\Br=X(\A_k)^\Br \times Y(\A_k)^\Br.$$}

\noindent The key topological fact behind our proof of Theorem C is this:
for any path-connected non-empty CW-complexes $X$ and $Y$, and 
any commutative ring $R$ with $1$ there is 
a canonical isomorphism
$$\H^2(X\times Y,R)=\H^2(X,R)\oplus\H^2(Y,R)\oplus 
\big(\H^1(X,R)\otimes_R\H^1(Y,R)\big).$$
See Proposition \ref{second} for this exercise in algebraic topology.
(This formula does not generalise to the third cohomology group, see
Remark \ref{counterex}.) The proof of Theorem C uses
Theorem \ref{kun0} that gives a similar result for the \'etale cohomology
of connected varieties over $\bar k$.

T. Schlank and Y. Harpaz, using \'etale homotopy of Artin and Mazur, 
recently proved a statement similar 
to our Theorem C where the Brauer--Manin set is replaced by
the {\it \'etale} Brauer--Manin set. In their result
the varieties
$X$ and $Y$ do not need to be proper, see \cite[Cor. 1.3]{SH}.

The first named author thanks the University of Tel Aviv, where a
part of this paper was written, for hospitality. He is grateful to
J.-L. Colliot-Th\'el\`ene, B. Kahn and A. P\'al for useful discussions.

%This should have consequences for the \'etale Brauer--Manin obstruction.

\section{Preliminaries}
\begin{sect}
{\bf Notation and conventions}. In this paper `almost all'
means `all but finitely many'.
If $B$ is an abelian group, we write
$B_\tors$ for the torsion subgroup of $B$. Let
$B/_\tors:=B/B_\tors$. If $\ell$  is a prime, then $B(\ell)$ is the
subgroup of $B_\tors$ consisting of the elements whose order is a
power of $\ell$, and $B({\non}\ell)$ is the subgroup of $B_\tors$
consisting of the elements whose order is {\it not} divisible by
$\ell$. If $m$ is a positive integer, then $B_m$ is the kernel of
the multiplication by $m$ in $B$.
\end{sect}

\begin{sect} 
{\bf Tate modules}. Let us recall some useful
elementary statements that are due to Tate \cite{tateB,Tate76}.
Let $B$ be an abelian group. The projective limit of the groups
$B_{\ell^n}$ (where the transition maps are the multiplications by
$\ell$) is called the $\ell$-adic Tate module of $B$, and is
denoted by $T_{\ell}(B)$. This limit carries a natural structure
of a  $\Z_{\ell}$-module; there is a natural injective map
$T_{\ell}(B)/\ell \hookrightarrow B_{\ell}$. One may easily check
that $T_{\ell}(B)_{\ell}=0$, and hence $T_{\ell}(B)$ is
torsion-free.

Let us assume that $B_{\ell}$ is finite. Then all the $B_{\ell^n}$
are obviously finite, and $T_{\ell}(B)$ is finitely generated by
Nakayama's lemma. Therefore, $T_{\ell}(B)$ is isomorphic to
$\Z_{\ell}^r$ for some non-negative integer $r\le
\dim_{\F_{\ell}}(B_{\ell})$. Moreover, $T_{\ell}(B)=0$ if and
only if $B(\ell)$ is finite. We denote by $V_{\ell}(B)$ the
$\Q_{\ell}$-vector space
$T_{\ell}(B)\otimes_{\Z_{\ell}}\Q_{\ell}$. Clearly,
$V_{\ell}(B)=0$ if and only if $B(\ell)$ is finite.

If $A$ is an abelian variety over a field $k$, and $\ell$ is a prime
different from $\fchar(k)$, $n=\ell^i$, then we write $A_n$ for the
kernel of multiplication by $n$ in $A(\bar{k})$. The group $A_n$
is a free $\Z/n$-module of rank $2\dim(A)$ equipped with the
natural structure of a $\Gamma$-module \cite{Mu}. We write
$T_{\ell}(A)$ for $T_{\ell}(A(\bar{k}))$, and $V_{\ell}(A)$ for
$V_{\ell}(A(\bar{k}))$, respectively. The $\Q_{\ell}$-vector space
$V_{\ell}(A)$ has dimension $2\dim(A)$ and carries the natural
structure of a $\Gamma$-module.
\end{sect}

\begin{sect} \label{1.3}
{\bf The Kummer sequence and the Picard variety}.
Let $X$ be a smooth, projective and geometrically integral
variety over $k$. 
For a positive integer $n$ coprime to $\fchar(k)$ we have
the Kummer exact sequence of sheaves of abelian groups 
in \'etale topology:
$$0\to\mu_n\to\G_m\to\G_m\to 0.$$
Recall that $\H_{\et}^1(\ov X,\G_m)=\Pic(\ov X)$.
Thus the Kummer sequence gives rise to an isomorphism of $\Ga$-modules
\begin{equation}
\H_{\et}^1(\ov X,\mu_n)=\Pic(\ov X)_n.\label{i1}
\end{equation}
Let $\Pic^0(\ov X)$ be the $\Gamma$-submodule
of $\Pic(\ov X)$ consisting of the classes of divisors algebraically
equivalent to $0$. By definition, the N\'eron--Severi group of $\ov X$
is the quotient $\Gamma$-module $\NS(\ov X)=\Pic(\ov X)/\Pic^0(\ov X)$.
The abelian group $\NS(\ov X)$ is finitely generated by a theorem of
N\'eron and Severi.

Let $A$ be the {\it Picard variety} of $X$, see \cite{Lang}.
Then $A$ is an abelian variety
over $k$ such that the $\Ga$-module $A(\bar k)$ is identified with
$\Pic^0(\ov X)$. 
Since the multiplication by $n$ is a surjective endomorphism of $A$,
we have an exact sequence of $\Gamma$-modules
\begin{equation}
0\to A_n\to \Pic(\ov X)_n \to \NS(\ov X)_n\to 0.\label{e3}
\end{equation}
Setting $n=\ell^m$, where $\ell\not=\fchar(k)$ is a prime,
we deduce from (\ref{i1}) and (\ref{e3})
a canonical isomorphism of $\Ga$-modules
\begin{equation}
\H_{\et}^1(\ov X,\Z_\ell(1))=T_\ell(A)=T_\ell(\Pic(\ov X)).\label{i2}
\end{equation}
In particular, this is a free $\Z_\ell$-module of finite rank.

Again, by surjectivity of the multiplication by $\ell^m$ on 
$A(\bar k)=\Pic^0(\ov X)$ we obtain from the Kummer sequence 
the following exact sequence of $\Ga$-modules
\begin{equation}
0\to\NS(\ov X)/\ell^m\to \H^2_\et(\ov X,\mu_{\ell^m})\to 
\Br(\ov X)_{\ell^m}\to 0.
\label{e5}
\end{equation}
Passing to the projective limit in $m$ gives rise to the
well known exact sequence
\begin{equation}
0\to\NS(\ov X)\otimes\Z_\ell\to \H^2_\et(\ov X,\Z_\ell(1))\to 
T_\ell(\Br(\ov X))\to 0.
\label{e6}
\end{equation}
It shows that the torsion subgroup of $\H^2_\et(\ov X,\Z_\ell(1))$
coincides with $\NS(\ov X)(\ell)$, cf. \cite[Sect. 2.2]{SZJAG}.
In particular, if $\ell$ does not divide the order of the torsion
subgroup of $\NS(\ov X)$, then $\H^2_\et(\ov X,\Z_\ell(1))$
is a free $\Z_\ell$-module of finite rank.

\end{sect}

\begin{sect} \label{1.4}
{\bf Products of varieties}.
Let $X$ and $Y$ be smooth, projective and geometrically integral
varieties over $k$.
We have the natural projection maps
$$\pi_X:X\times Y\to X, \quad \pi_Y:X\times Y\to Y.$$
We denote by the same symbols
the projections $\ov X\times \ov Y\to \ov X$ and
$\ov X\times \ov Y\to \ov Y$. Fixing $\bar k$-points
${x}_0\in X(\bar{k})$ and ${y}_0\in Y(\bar{k})$,
we define closed embeddings
$$q_{{y}_0}: \ov X=\ov X\times  {y}_0\hookrightarrow
\ov X\times \ov Y, \quad  q_{{x}_0}: \ov Y=x_0\times \ov Y
\hookrightarrow  \ov X\times  \ov Y.$$
Then $\pi_X q_{{y}_0}=\id_{\ov X}$ and
$\pi_Y q_{{x}_0}=\id_{\ov Y}$.
On the other hand,
 $$\pi_Y q_{{y}_0}(\ov X)=y_0 \subset \ov Y, \quad
 \pi_X q_{{x}_0}(\ov Y)=x_0 \subset  \ov X.$$

Let $\sF$ be an \'etale sheaf defined by a commutative
$k$-group scheme, see \cite[Cor. II.1.7]{EC}. For example,
$\sF$ can be the sheaf defined by the multiplicative group $\G_m$,
or by the finite $k$-groups $\Z/n$ or $\mu_n$, where
$n$ is not divisible by the characteristic of $k$.
The induced map
$\pi_X^{*}: \H_{\et}^i(\ov X,\sF)\to \H_{\et}^i(\ov X\times\ov Y,\sF)$
is a homomorphism of $\Gamma$-modules, whereas
$q_{{x}_0}^{*}:\H_{\et}^i(\ov X\times\ov Y,\sF) \to \H_{\et}^i(\ov Y,\sF)$
is {\it a priori} only a homomorphism of abelian groups.
If $x_0\in X(k)$, then $q_{{x}_0}^{*}$ is also a homomorphism of $\Ga$-modules.

The next proposition easily follows from definitions and the above
considerations.
\end{sect}

\begin{prop}
\label{obvious} For any $i\geq 1$ we have the following statements.

{\rm (i)} The induced maps
$$\pi_X^{*}: \H_{\et}^i(\ov X,\sF)\to \H_{\et}^i(\ov X\times\ov Y,\sF),
\
\pi_Y^{*}: \H_{\et}^i(\ov Y,\sF)\to
\H_{\et}^i(\ov X\times\ov Y,\sF)$$
are injective homomorphisms of $\Gamma$-modules.

{\rm (ii)} The induced maps $q_{{y}_0}^{*}$ and $q_{{x}_0}^{*}$
define isomorphisms of abelian groups
$$q_{{y}_0}^{*}: \pi_X^{*}(\H_{\et}^i(\ov X,\sF))\to
\H_{\et}^i(\ov X,\sF), \ q_{{x}_0}^{*}:
\pi_Y^{*}(\H_{\et}^i(\ov Y,\sF))\to \H_{\et}^i(\ov Y,\sF).$$

{\rm (iii)} The subgroup $\pi_X^{*}(\H_{\et}^i(\ov X,\sF))$
lies in the kernel of
$$q_{{x}_0}^{*}:\H_{\et}^i(\ov X\times\ov Y,\sF) \to \H_{\et}^i(\ov Y,\sF),$$
and, similarly, $\pi_Y^{*}(\H^i(\ov Y,\sF))$ lies in the
kernel of
$$q_{{y}_0}^{*}:\H_{\et}^i(\ov X\times\ov Y,\sF) \to \H_{\et}^i(\ov X,\sF).$$
Hence we have $\pi_X^{*}(\H_{\et}^i(\ov X,\sF))\cap
\pi_Y^{*}(\H_{\et}^i(\ov Y,\sF))=0$.

{\rm (iv)} The map $(a,b) \mapsto \pi_X^{*}(a)+\pi_Y^{*}(b)$
defines an injective homomorphism of $\Ga$-modules
$$\H_{\et}^i(\ov X,\sF)\oplus \H_{\et}^i(\ov Y,\sF)\to
\H_{\et}^i(\ov X\times\ov Y,\sF).$$
\end{prop}

\begin{sect}{\bf Picard groups of products of varieties.}
We identify $\Pic(\ov X)\oplus \Pic(\ov Y)$ with its image in
$\Pic(\ov X\times \ov Y)$. It is well known that
$$\Pic^0(\ov X\times \ov Y)=\Pic^0(\ov X)\oplus \Pic^0(\ov Y).$$
%reference
Let $A$ be the Picard variety of $X$, and let $B$
be the Picard variety of $Y$.
The dual abelian variety $A^t$ of $A$ is the {\it Albanese
variety} of $X$. When $X(k)\not=\emptyset$, 
the choice of a point $x_0\in X(k)$
defines a morphism ${\rm Alb}_{x_0}:X\to A^t$ that sends $x_0$ to $0$.
The pair $(A^t,{\rm Alb}_{x_0})$ can be characterized by the universal property
that any morphism from $X$ to an abelian variety $A'$
that sends $x_0$ to $0$ is the composition of ${\rm Alb}_{x_0}$ and
a morphism of abelian varieties $A^t\to A'$. See  \cite{Mu}, \cite{Lang}
for more details. It is clear that the Albanese variety of $\ov X$
is $\ov{A^t}$.

\begin{prop} \label{diag}
We have a commutative diagram of $\Gamma$-modules with
exact rows and columns, where the exact sequence in the bottom row is split:
$$\begin{array}{ccccccccc}
&&0&&0&&&&\\
&&\downarrow&&\downarrow&&&&\\
&&A(\bar k)\oplus B(\kbar)&=&A(\bar k)\oplus B(\kbar)&&&&\\
&&\downarrow&&\downarrow&&&&\\
0&\to& \Pic(\ov X)\oplus \Pic(\ov Y)&\to& \Pic(\ov X\times \ov Y)&\to&
\Hom(\ov{B^t},\ov A)&\to&0\\
&&\downarrow&&\downarrow&&||&&\\
0&\to& \NS(\ov X)\oplus \NS(\ov Y)&\to& \NS(\ov X\times \ov Y)&\to&
\Hom(\ov{B^t},\ov A)&\to&0\\
&&\downarrow&&\downarrow&&&&\\
&&0&&0&&&&
\end{array}$$
If $(X\times Y)(k)\not=\emptyset$, then the exact sequence in the
middle row is also split.
\end{prop}
\noindent{\bf Proof.} Choose a $\bar k$-point $(x_0,y_0)$ in $X\times Y$, and
let $P_{x_0,y_0}$ be the kernel of the group homomorphism
$$\Pic(\ov X\times \ov Y)\to \Pic(\ov X)\oplus \Pic(\ov Y), \quad
L \mapsto (q_{{y}_0}^{*}L,q_{{x}_0}^{*}L).$$
Let $N_{x_0,y_0}$ be the image of $P_{x_0,y_0}$ in $\NS(\ov X\times\ov Y)$.
By Proposition \ref{obvious} the intersection of
$P_{x_0,y_0}$ with $\Pic(\ov X)\oplus \Pic(\ov Y)$ inside
$\Pic(\ov X\times \ov Y)$ is zero, hence the natural surjective
map $P_{x_0,y_0}\to N_{x_0,y_0}$ is an isomorphism of abelian groups.

For any $L\in P_{x_0,y_0}$ we have
$q_{{y}_0}^{*}L=0$, hence $q_{y}^{*}L\in \Pic^0(\ov X)$
for any $y\in Y(\bar k)$. Thus $N_{x_0,y_0}$ is the kernel of the group
homomorphism
$$\NS(\ov X\times \ov Y)\to \NS(\ov X\times y)\oplus \NS(x\times\ov Y)$$
for any $x\in X(\bar k)$ and $y\in Y(\bar k)$. In particular, $N_{x_0,y_0}$
does not depend on the choice of $(x_0,y_0)$, so we can
drop $x_0$ and $y_0$, and write $N=N_{x_0,y_0}$.
It follows that $N$ is a Galois submodule of
$\NS(\ov X\times \ov Y)$, so that we have a decomposition
of $\Gamma$-modules
$$\NS(\ov X\times \ov Y)=\NS(\ov X)\oplus \NS(\ov Y)\oplus N.$$
It remains to show that $\Gamma$-modules $N$ and
$\Hom(\ov{B^t},\ov A)$ are canonically isomorphic.

The Poincar\'e sheaf ${\mathcal P}_X$ on
$A^t\times A$ is a certain canonical invertible sheaf that restricts
trivially to both $\{0\}\times A$ and $A^t\times \{0\}$, see \cite{Mu}, \cite{Lang}.
Every morphism of abelian varieties $u:\ov B^t\to\ov A$
gives rise to the invertible sheaf
$({\rm Alb}_{x_0},u\,{\rm Alb}_{y_0})^{*}({\mathcal P}_X)$ on 
$\ov X\times\ov Y$, whose isomorphism class is in $P_{x_0,y_0}$.
It is well known that this defines a group isomorphism
\begin{equation}
\Hom(\ov B^t,\ov A)\tilde\lra P_{x_0,y_0}.\label{dc}
\end{equation}
%we need a reference for this fact!
The Poincar\'e sheaf is defined over $k$ so from (\ref{dc}) we deduce
a canonical isomorphism of $\Gamma$-modules
$\Hom(\ov B^t,\ov A)\tilde\lra N$. The last statement of the proposition
is clear: it is enough to choose $(x_0,y_0)\in (X\times Y)(k)$.
QED

\medskip

The following corollary is well known. See Proposition \ref{second} below
for a topological analogue.

\begin{cor} \label{ku1}
Let $n$ be a positive integer not divisible by $\fchar(k)$.
Then we have a canonical decomposition of $\Ga$-modules
\begin{equation}
\H_{\et}^1(\ov X\times \ov Y,\mu_n)=\H_{\et}^1(\ov X,\mu_n)\oplus
\H_{\et}^1(\ov Y,\mu_n).\label{Ku1}
\end{equation}
\end{cor}
\noindent{\bf Proof}. The middle row of the diagram of Proposition \ref{diag}
gives an isomorphism of $\Ga$-modules $\Pic(\ov X)_n\oplus\Pic(\ov Y)_n$
and $\Pic(\ov X\times\ov Y)_n$. It remains to use the canonical isomorphism
(\ref{i1}). QED

\begin{rem} \label{1.9}
The abelian group $\Hom(\ov{B^t},\ov A)$ is
finitely generated and torsion-free, hence
$\H^1(k,\Hom(\ov{B^t},\ov A))$ is finite. It follows that
the cokernel of the natural map
$$\H^1(k,\Pic(\ov X))\oplus \H^1(k,\Pic(\ov Y)) \lra
\H^1(k,\Pic(\ov X\times \ov Y)),$$
and the kernel of the natural map
$$\H^2(k,\Pic(\ov X))\oplus\H^2(k,\Pic(\ov Y)) \lra
\H^2(k,\Pic(\ov X\times \ov Y))$$
are both finite.
\end{rem}

\end{sect}

\section{K\"unneth decompositions} \label{K}

\begin{sect}{\bf K\"unneth decomposition with coefficients in a field}.
We continue to assume that $X$ and $Y$ are smooth, projective and 
geometrically integral varieties over $k$.
Let $\ell\not=\fchar(k)$ be a prime.
We have the K\"unneth decomposition of $\Gamma$-modules
$$\H_{\et}^2(\ov X\times
\ov Y,\Q_{\ell})=\H_{\et}^2(\ov X,\Q_{\ell})\oplus
\H_{\et}^2(\ov Y,\Q_{\ell})\oplus
\big(\H_{\et}^1(\ov Y,\Q_{\ell})\otimes_{\Q_{\ell}}
\H_{\et}^1(\ov X,\Q_{\ell})\big),$$
see \cite[Cor. VI.8.13]{EC}. 
From (\ref{i2}) we have a canonical isomorphism
$$\H_{\et}^1(\ov X,\Q_{\ell}(1))=V_{\ell}(A).$$
When $n$ is a positive integer coprime to $\fchar(k)$, 
the non-degeneracy of the Weil pairing gives rise to a canonical 
isomorphism of Galois modules
$B_n=\Hom(B^t_n,\mu_n)$,
and hence to a canonical isomorphism
$$\H_{\et}^1(\ov Y,\Q_{\ell})=V_\ell(B)(-1)=\Hom_{\Q_{\ell}}(V_{\ell}(B^t),\Q_{\ell}).$$
Therefore we have an isomorphism of $\Ga$-modules \cite[p. 143]{tate}
$$\H_{\et}^1(\ov Y,\Q_{\ell})\otimes_{\Q_{\ell}}
\H_{\et}^1(\ov X,\Q_{\ell}(1))=\Hom_{\Q_{\ell}}(V_{\ell}(B^t),V_{\ell}(A)),$$
and hence a decomposition of Galois modules
\begin{equation}
\H_{\et}^2(\ov X\times\ov Y,\Q_{\ell}(1))=\H_{\et}^2(\ov X,\Q_{\ell}(1))\oplus
\H_{\et}^2(\ov Y,\Q_{\ell}(1))\oplus
\Hom_{\Q_{\ell}}(V_{\ell}(B^t),V_{\ell}(A)).\label{d1}
\end{equation}

Our next result, Theorem \ref{kun0}, is probably well known to experts;
we give a proof as we could not find it in the literature.
As a motivation and for the sake of completeness we present 
a similar result for CW-complexes (we shall not need it
in the rest of the paper).

\end{sect}

\begin{prop} \label{second}
Let $X$ and $Y$ be non-empty path-connected CW-complexes.
For any commutative ring $R$ with $1$
we have canonical isomorphisms of abelian groups
$$\H^1(X\times Y,R)=\H^1(X,R)\oplus\H^1(Y,R)$$
and
$$\H^2(X\times Y,R)=\H^2(X,R)\oplus\H^2(Y,R)\oplus 
\big(\H^1(X,R)\otimes_R\H^1(Y,R)\big).$$
\end{prop}
\noindent{\bf Proof.} To simplify notation we write
$\H_n(X)$ for $\H_n(X,\Z)$.
The universal coefficients theorem \cite[Thm. 3.2]{Hatch}
gives the following (split) exact sequence of abelian groups
\begin{equation}
0\to \Ext(\H_{n-1}(X),R)\to \H^n(X,R)\to
\Hom(\H_n(X),R)\to 0.\label{univ}
\end{equation}
Since $X$ is non-empty and path-connected we have 
$\H_0(X)=\Z$, see \cite[Prop. 2.7]{Hatch}. This 
gives a canonical isomorphism
\begin{equation}
\H^1(X,R)=\Hom(\H_1(X),R).\label{h_1}
\end{equation} 
The K\"unneth formula for homology is 
$$
0\to \bigoplus_{i=0}^n \big(\H_i(X)\otimes \H_{n-i}(Y)\big)
\to \H_n(X\times Y)\to 
\bigoplus_{i=0}^{n-1}\Tor(\H_i(X),\H_{n-1-i}(Y))\to 0,
$$
see \cite[Thm. 3.B.6]{Hatch}. We deduce from it
canonical isomorphisms
\begin{equation}
\H_1(X\times Y)=\H_1(X)\oplus \H_1(Y)
\label{h1}
\end{equation}
and 
\begin{equation}
\H_2(X\times Y)=\H_2(X)\oplus\H_2(Y)\oplus 
\big(\H_1(X)\otimes\H_1(Y)\big).
\label{h2}
\end{equation}
Our first isomorphism follows from (\ref{h_1}) and (\ref{h1}).

The exact sequence (\ref{univ}) for $n=2$ gives rise to 
the following commutative diagram
$$\begin{array}{ccccccccc}
0&\to& \Ext(\H_1(X)\oplus\H_1(Y),R)&\to& \H^2(X,R)\oplus \H^2(Y,R)&\to&
\Hom(\H_2(X)\oplus\H_2(Y),R)&\to& 0\\
&&\uparrow&&\uparrow&&\uparrow&&\\
0&\to& \Ext(\H_1(X\times Y),R)&\to& \H^2(X\times Y,R)&\to&
\Hom(\H_2(X\times Y),R)&\to& 0\end{array}$$
By (\ref{h1}) the left vertical arrow is an isomorphism,
hence the kernels of the other two vertical arrows
are isomorphic. Hence, by (\ref{h2}), the kernel of the middle
vertical map is isomorphic to
$\Hom(\H_1(X)\otimes\H_1(Y),R)$, which by (\ref{h_1})
is isomorphic to $\H^1(X,R)\otimes_R\H^1(Y,R)$. 
Moreover, $\H^2(X,R)$ and $\H^2(Y,R)$ are direct factors of 
$\H^2(X\times Y,R)$, so our second isomorphism follows. QED

\begin{rem} \label{counterex}
Let $X=\R\P^2$. Then $\H_1(X)=\Z/2$ and $\H_n(X)=0$ for $n\geq 2$.
From the universal coefficients theorem (\ref{univ})
we obtain $\H^1(X,\Z)=0$, $\H^2(X,\Z)=\Z/2$ and $\H^n(X,\Z)=0$
for $n\geq 3$, cf. \cite[Ex. 3.9]{Hatch}. Combining the
calculation of homology of $X^2$ in \cite[Ex. 3.B.3]{Hatch} 
with the universal coefficients theorem we obtain 
$$\H^3(X^2,\Z)=\Z/2\ \not=\ \bigoplus_{i=0}^3 \big(\H^i(X,\Z)\otimes
\H^{3-i}(X,\Z)\big).$$
This shows that Proposition \ref{second} does not generalise
to the third cohomology group, at least when $R=\Z$.
\end{rem}

\begin{sect}{\bf The type of a torsor}. \label{ttype}
After this digression into algebraic topology we return to 
smooth, projective and geometrically integral
varieties $X$ and $Y$ over a field $k$. We now
introduce some notation. 
Let $S_X$ be the finite commutative $k$-group of multiplicative type
whose Cartier dual $\hat S_X:=\Hom_{\bar k{\rm-gr.}}(S_X,\G_m)$ is
$$\hat S_X=\H_\et^1(\ov X,\mu_n)=\Pic(\ov X)_n.$$ 
Then we have a canonical identification
\begin{equation}
\Hom(S_X,\Z/n)=\H_\et^1(\ov X,\Z/n). \label{iden}
\end{equation}
For any finite $\bar k$-group scheme $G$ of multiplicative type
annihilated by $n$ 
we have a canonical isomorphism, functorial in $X$ and $G$:
$$\tau_G:\H_\et^1(\ov X,G)\tilde\lra\Hom(\hat G,\Pic(\ov X))=
\Hom(\hat G,\hat S_X),$$
see \cite[Cor. 2.3.9]{S} (this uses the assumption that
$\ov X$ is projective and connected, and so has no non-constant invertible
regular functions).
It can be defined via the natural pairing 
\begin{equation}
\H^1_\et(\ov X,G) \times \hat G \to \H^1_\et(\ov X,\mu_n)=\hat S_X, \label{17}
\end{equation}
see \cite[Section 2.3]{S}.
If $Z/\ov X$ is a torsor under $G$, then the associated homomorphism 
$\tau_G(Z):\hat G\to \hat S_X$ is called the {\it type} of $Z/\ov X$. 
If we take $G=S_X$,
then there exists a torsor $\T_{\ov X}/\ov X$ under $S_X$, unique
up to isomorphism, whose type is the identity map. Thus there is
a well defined class $[\T_{\ov X}]\in \H^1_\et(\ov X,S_X)$. 
This class can be
used to describe $\tau_G^{-1}$ explicitly. For 
$\varphi\in \Hom(\hat G,\hat S_X)$ let
$\hat\varphi\in\Hom(S_X,G)$ be the homomorphism that corresponds to
$\varphi$ under the identification
$$\Hom(\hat G,\hat S_X)=\Hom(S_X,G).$$
The functoriality of $\tau_G$ in $G$ implies that 
$\tau_G^{-1}(\varphi)$ is the push-forward $\hat\varphi_*[\T_{\ov X}]$,
which can also be defined as the class of the $\ov X$-torsor
$(\T_{\ov X}\times_k G)/S_X$.

If we take $G=S_X$ in (\ref{17}) we obtain a natural pairing
$$\H^1_\et(\ov X,S_X) \times \hat S_X \to \H^1_\et(\ov X,\mu_n)=\hat S_X.$$
The definition of $\T_{\ov X}$ implies that pairing with
the class $[\T_{\ov X}]$ gives the identity map on $\hat S_X$.
After a twist we obtain a natural pairing 
$$\H^1_\et(\ov X,S_X) \times \H^1_\et(\ov X,\Z/n) \to \H^1_\et(\ov X,\Z/n);$$
moreover, pairing with $[\T_{\ov X}]$ gives the identity map on
$\H^1_\et(\ov X,\Z/n)$.

\end{sect}

\begin{rem} \label{rem14}
There is a natural cup-product map
$$\H^1_\et(\ov X,S_X)\otimes\H^1_\et(\ov Y,S_Y)\lra
\H^2_\et(\ov X\times \ov Y,S_X\otimes S_Y).$$
Let us denote the image of $[\T_{\ov X}]\otimes[\T_{\ov Y}]$ by 
$[\T_{\ov X}]\cup[\T_{\ov Y}]$.
From (\ref{iden}) we obtain a natural pairing
$$\H^2_\et(\ov X\times \ov Y,S_X\otimes S_Y)\ \times\ 
\H_\et^1(\ov X,\Z/n)\otimes\H_\et^1(\ov Y,\Z/n)
\ \lra\ \H^2_\et(\ov X\times \ov Y,\Z/n).$$
Since pairing with $[\T_{\ov X}]$ induces identity on $\H_\et^1(\ov X,\Z/n)$,
and similarly for $Y$, we see that pairing with 
$[\T_{\ov X}]\cup[\T_{\ov Y}]$ gives the cup-product map 
$$\cup:\ \H_\et^1(\ov X,\Z/n)\otimes\H_\et^1(\ov Y,\Z/n)\to 
\H^2_\et(\ov X\times \ov Y,\Z/n).$$
\end{rem}

\begin{thm} \label{kun0}
Let $n$ be a positive integer coprime to $\fchar(k)$.
Then the homomorphism of $\Ga$-modules
$$
\H_{\et}^2(\ov X,\Z/n)\oplus
\H_{\et}^2(\ov Y,\Z/n)\oplus 
\big(\H_{\et}^1(\ov X,\Z/n)\otimes\H_{\et}^1(\ov Y,\Z/n)\big)
\lra \H_{\et}^2(\ov X\times \ov Y,\Z/n)
$$
given by $\pi_X^*$ on the first factor, by $\pi_Y^*$ on the second factor,
and by the cup-product on the third factor,
is an isomorphism.
\end{thm}
\noindent{\bf Proof.} 
It is enough to establish this decomposition at the level of abelian groups.
Choose $x_0\in X(\bar k)$, $y_0\in Y(\bar k)$.
Using the notation of Section \ref{1.4} we define
$$\H^2_\et(\ov X\times \ov Y,\Z/n)_{\rm prim}=
\Ker[(q_{y_0}^*,q_{x_0}^*): \H^2_\et(\ov X\times \ov Y,\Z/n)\to
\H^2_\et(\ov X,\Z/n)\oplus\H^2_\et(\ov Y,\Z/n)].$$

The \'etale (or Zariski) sheaf $R^q\pi_{X*}(\Z/n)$ is the constant sheaf
associated with the finite abelian group $\H_\et^q(\ov Y,\Z/n)$
(for example, by the proper base change theorem \cite[Cor. VI.2.3]{EC}).
Thus we have the following Leray spectral sequence 
\begin{equation}
E^{p,q}_2=\H_\et^p(\ov X,\H_\et^q(\ov Y,\Z/n))\Rightarrow 
\H^{p+q}_\et(\ov X\times \ov Y,\Z/n).
\label{a0}
\end{equation}
We have seen in Proposition \ref{obvious} 
that the maps $\pi_X^*$ and $\pi_Y^*$ make
the abelian groups
$\H^m_\et(\ov X,\Z/n)$ and $\H^m_\et(\ov Y,\Z/n)$ direct
factors of $\H^m_\et(\ov X\times \ov Y,\Z/n)$, for all $m\geq 1$.
By the standard theory
of spectral sequences this gives a canonical isomorphism
$$\beta:\H^2_\et(\ov X\times \ov Y,\Z/n)_{\rm prim} 
\ \tilde\lra\ \H_\et^1(\ov X,\H_\et^1(\ov Y,\Z/n)).$$

Taking $G=\hat S_Y$ in (\ref{17}) we get an isomorphism
$\tau_{\hat S_Y}:\H_\et^1(\ov X,\hat S_Y)\tilde\lra\Hom(S_Y,\hat S_X)$. 
Using (\ref{iden}), after a twist we obtain an isomorphism
$$
\tau:\H_\et^1(\ov X,\H_\et^1(\ov Y,\Z/n))\ \tilde\lra\ 
\H_\et^1(\ov X,\Z/n)\otimes\H_\et^1(\ov Y,\Z/n).
$$
This gives some isomorphism as in the statement
of the theorem.
To complete the proof we need to check that for any 
$x\in \H_\et^1(\ov X,\Z/n)$
and any $y\in \H_\et^1(\ov Y,\Z/n)$ we have
$$x \cup y=\beta^{-1}\tau^{-1}(x\otimes y).$$

We have seen above that $\tau^{-1}(x\otimes y)$
is the push-forward of $[\T_{\ov X}]$ by the map
$S_X\to \Hom(S_Y,\Z/n)$ defined by
$$x\otimes y\in \H^1_\et(\ov X,\Z/n)\otimes\H^1_\et(\ov Y,\Z/n)=
\Hom(S_X\otimes S_Y,\Z/n).$$
In other words, $\tau^{-1}(x\otimes y)$ is obtained by
pairing $[\T_{\ov X}]$ with $x\otimes y$. In view of Remark
\ref{rem14} in order to finish the proof, it remains to check that 
$\beta^{-1}$ can be described via the pairing
$$\H_\et^1(\ov X,\Hom(S_Y,\Z/n))\times \H^1_\et(\ov Y,S_Y)\lra 
\H^2_\et(\ov X\times \ov Y,\Z/n),$$
namely, as pairing with $[\T_{\ov Y}]$. This calculation is more or less
standard, cf. \cite[Thm. 4.1.1]{S} or, more recently, \cite[Thm. 1.4]{HS}.

Let us write $\D(Z)$ for the bounded derived
category of \'etale sheaves of abelian groups on a variety $Z$. Let 
$\bR{\pi_X}_*:\D(\ov X\times \ov Y)\to \D(\ov X)$
be the derived functor of $\pi_X$. Let 
$\rho:\ov Y\to \Spec(\bar k)$ be the structure
morphism, and let 
$\bR \rho_*:\D(\ov Y)\to \D({\rm Ab})$
be the corresponding derived functor to 
the bounded derived category of the category
of abelian groups ${\rm Ab}$. Each of these derived categories has
the canonical truncation functors $\tau_{\leq m}$.
We need to recall the definition of the type map
of a group $G$ of multiplicative type, see
\cite[Section 2.3]{S}. This is the composed map
\begin{equation}
\H^1_\et(\ov Y,G)\to 
\Ext^1(\hat G,\tau_{\leq 1}\bR \rho_*\G_m)\to \Hom(\hat G,\Pic(\ov Y)).
\label{type}
\end{equation}
The Hom- and Ext-groups without subscript are taken in ${\rm Ab}$
or $\D({\rm Ab})$.
The second map in (\ref{type}) is induced by the obvious exact triangle 
in $\D({\rm Ab})$
$$\bar k^*\to\tau_{\leq 1}\bR \rho_*\G_m\to (\Pic \ov Y)[-1],$$ 
where we used the facts that $\H^0_\et(\ov Y,\G_m)=\bar k^*$, since
$\ov Y$ is reduced and connected, and $\H^1_\et(\ov Y,\G_m)=\Pic(\ov Y)$.
To define the first map in (\ref{type}) consider the local-to-global 
spectral sequence of Ext-groups 
$$E^{p,q}_2=\H^p_\et(\ov Y,\mathcal{E}xt^q_{\ov Y}(\hat G,\G_m))\Rightarrow
\Ext^{p+q}_{\ov Y}(\hat G,\G_m).$$
It completely degenerates since $\mathcal{E}xt^q_{\ov Y}(\hat G,\G_m)=0$
for $q\geq 1$, thus giving an isomorphism $\H^q_\et(\ov Y,G)\tilde\lra
\Ext^q_{\ov Y}(\hat G,\G_m)$
\cite[Lemma 2.3.7]{S}. It remains to use the identities
$$\Ext^q_{\ov Y}(\hat G,\G_m)=\Ext^q(\hat G,\bR\rho_*\G_m)=
\Ext^q(\hat G,\tau_{\leq q}\bR \rho_*\G_m)$$
stemming from the fact that
$\bR\Hom_{\ov Y}(\rho^*\hat G, \cdot)=\bR\Hom(\hat G,\bR\rho_*(\cdot))$.
When $G$ is annihilated by $n$, the image of the type map lies in
$\Hom(\hat G,(\Pic \ov Y)_n)$, and thus $\tau_G$ can be written as the composition
of the maps
$$\H^1_\et(\ov Y,G)\to 
\Ext^1(\hat G,\tau_{\leq 1}\bR \rho_*\mu_n)\to \Hom(\hat G,(\Pic \ov Y)_n).$$
We claim that these maps fit into the following commutative diagram of pairings:
$$\begin{array}{ccccc}
\H^1_\et(\ov X,\hat G)&\times&\H^1_\et(\ov Y,G)&\to&
\H^2_\et(\ov X\times \ov Y,\mu_n)\\
||&&\downarrow&&\uparrow\\
\H^1_\et(\ov X,\hat G)&\times&\Ext^1(\hat G,\tau_{\leq 1}\bR \rho_*\mu_n)
&\to&\H^2_\et(\ov X,\tau_{\leq 1}\bR {\pi_X}_*\mu_n)\\
||&&\downarrow&&\downarrow\\
\H^1_\et(\ov X,\hat G)&\times&\Hom(\hat G,(\Pic \ov Y)_n)
&\to&\H^1_\et(\ov X,(\Pic \ov Y)_n)
\end{array}$$
The first pairing is induced by the obvious pairing
$$\pi_X^*\hat G\ \times \ 
{\mathcal H}om_{\ov Y}(\pi_Y^*\hat G,\mu_{n,\ov Y})
\ \lra \ \mu_{n,\ov X\times\ov Y},$$
by taking cohomology of $\ov X$, $\ov Y$, $\ov X\times\ov Y$, 
respectively. If instead of 
$R^1\rho_*({\mathcal H}om_{\ov Y}(\cdot,\cdot))$
we consider $R^1(\rho_*{\mathcal H}om_{\ov Y})(\cdot,\cdot)$
we get the second pairing. This explains the compatibility of
the two upper pairings. The third pairing
is the natural one, so that the compatibility of 
the two lower pairings is clear.

Now take $G=S_Y$, so that $\hat G=(\Pic \ov Y)_n$. By pairing with 
$[\T_{\ov Y}]$,
after a twist we obtain a map 
$$\gamma:\H^1_\et(\ov X,\H^1(\ov Y,\Z/n))\to \H^2_\et(\ov X\times \ov Y,\Z/n),$$
which factors through the injective map 
$$\H^2_\et(\ov X,\tau_{\leq 1}\bR {\pi_X}_*\mu_n) \to 
\H^2_\et(\ov X\times \ov Y,\Z/n).$$
Since $\bar k$ is separably closed we have
$\H^1_\et(y_0,G)=\H^1_\et(\bar k,G)=0$, so that $q_{y_0}^*\gamma=0$.
A similar argument gives $q_{x_0}^*\gamma=0$, thus
$\Im(\gamma)\subset\H^2_\et(\ov X\times \ov Y,\Z/n)_{\rm prim}$.
By the standard theory of spectral sequences
the map $\beta$ is obtained from the right hand downward map
in the diagram
(after a twist). Since the type of $\T_{\ov Y}$ is the identity in 
$\Hom((\Pic \ov Y)_n,(\Pic \ov Y)_n)$, 
the commutativity of the diagram implies that $\beta\gamma=\id$. QED

\begin{cor} \label{kun}
Let $n$ be a positive integer coprime to $\fchar(k)$, $|\NS(\ov X)_{\rm tors}|$
and $|\NS(\ov Y)_{\rm tors}|$.
Then we have a canonical decomposition of $\Ga$-modules
\begin{equation}
\H_{\et}^2(\ov X\times \ov Y,\mu_n)=\H_{\et}^2(\ov X,\mu_n)\oplus
\H_{\et}^2(\ov Y,\mu_n)\oplus \Hom(B^t_n,A_n).\label{Ku}
\end{equation}
\end{cor}
\noindent{\bf Proof.} For any prime $\ell$ dividing $n$ we
have $\NS(\ov X)(\ell)=0$ and $\NS(\ov Y)(\ell)=0$. Thus 
from the isomorphism (\ref{i1}) and the exact sequence (\ref{e3})
we obtain canonical isomorphisms 
$$\H_{\et}^1(\ov X,\mu_{\ell^m})=A_{\ell^m}, \quad
\H_{\et}^1(\ov Y,\mu_{\ell^m})=B_{\ell^m},$$ 
for any $m\geq 1$.
From the non-degeneracy of the Weil 
pairing we deduce a canonical isomorphism of $\Gamma$-modules
$$\H^1(\ov Y,\Z/{\ell^m})=B_{\ell^m}(-1)=\Hom(B^t_{\ell^m},\Z/{\ell^m}).$$
We conclude that the Galois modules
$\H_{\et}^1(\ov Y,\Z/{\ell^m})\otimes\H_{\et}^1(\ov X,\mu_{\ell^m})$ 
and $\Hom(B^t_{\ell^m},A_{\ell^m})$ are canonically isomorphic. 
Hence, after a twist by $\mu_n$, the isomorphism of Theorem \ref{kun0} 
can be written as (\ref{Ku}). QED

\begin{sect} \label{1.12}
{\bf  First Chern classes}. Let $\ell$ be a prime different from $\fchar(k)$.
Tensoring (\ref{e6}) with $\Q_\ell$ we 
obtain the following exact sequence of $\Ga$-modules
\begin{equation}
0 \to \NS(\ov Z)\otimes\Q_{\ell} \to \H_{\et}^2(\ov Z,\Q_{\ell}(1)) \to
V_{\ell}(\Br(\ov Z)) \to 0, \label{odin}
\end{equation}
The injective maps from (\ref{odin}) and (\ref{e5}) are both called 
{\it the first Chern class} maps (see, e.g., \cite{EC}, VI.9):
$$c_1:\NS(\ov X)\otimes\Q_{\ell} \hookrightarrow \H_{\et}^2(\ov X,\Q_{\ell}(1)),
\quad \bar{c}_1:\NS(\ov X)/\ell^m \hookrightarrow \H_{\et}^2(\ov X,\mu_{\ell^m}).$$
Proposition \ref{diag} gives a natural isomorphism of Galois modules
$$\NS(\ov X\times \ov Y)=\NS(\ov X)\oplus
\NS(\ov Y)\oplus \Hom(\ov B^t,\ov A).$$ 
Since the maps $c_1$ and $\bar c_1$ are functorial in $X$, we see that
the map
$$c_1:\NS(\ov X\times \ov Y)\otimes \Q_{\ell} \hookrightarrow
\H_{\et}^2(\ov X\times \ov Y,\Q_{\ell}(1))$$
forms obvious commutative diagrams with the maps
$$c_1: \NS(\ov X)\otimes\Q_{\ell}\hookrightarrow
\H_{\et}^2(\ov X,\Q_{\ell}(1)), \quad
c_1: \NS(\ov Y)\otimes\Q_{\ell}\hookrightarrow \H_{\et}^2(\ov Y,\Q_{\ell}(1)),$$
$$c_1: \Hom(\ov B^t,\ov A)\otimes \Q_{\ell}\hookrightarrow
\Hom_{\Q_{\ell}}(V_{\ell}(B^t),V_{\ell}(A)).$$
Similarly, for $\ell$ coprime to $\fchar(k)$, $|\NS(\ov X)_{\rm tors}|$
and $|\NS(\ov Y)_{\rm tors}|$, the map
$$\bar{c}_1: \NS(\ov X\times \ov Y)/{\ell^m}\hookrightarrow
\H_{\et}^2(\ov X\times \ov Y,\mu_{\ell^m})$$
forms similar commutative diagrams with the maps
$$\bar{c}_1:\NS(\ov X)/{\ell^m}\hookrightarrow\H_{\et}^2(\ov X,\mu_{\ell^m}),
\quad
\bar{c}_1:\NS(\ov Y)/{\ell^m}\hookrightarrow   \H_{\et}^2(\ov Y,\mu_{\ell^m}),$$
$$\bar{c}_1:\Hom(\ov B^t,\ov A)/{\ell^m}\hookrightarrow \Hom(B^t_{\ell^m},A_{\ell^m})=
\H_{\et}^1(\ov Y,\Z/{\ell^m})\otimes\H_{\et}^1(\ov X,\mu_{\ell^m}).$$

\end{sect}

\begin{sect}
\label{BrauerAV}{\bf Brauer groups of products of varieties}.
Let $\ell$ be a prime different from $\fchar(k)$.
Applying (\ref{odin}) to $X$, $Y$ and $X \times Y$, using (\ref{d1}) 
and the compatibilities from Section \ref{1.12},
we obtain the following decomposition of Galois modules:
\begin{equation}V_{\ell}(\Br(\ov X\times\ov Y))=V_{\ell}(\Br(\ov X))\oplus
V_{\ell}(\Br(\ov Y)) \oplus
\big(\Hom_{\Q_{\ell}}(V_{\ell}(B^t),V_{\ell}(A))/\Hom(\ov B^t,\ov A)\otimes\Q_{\ell}\big).
\label{e1}\end{equation}
When $\ell$ is also coprime to $|\NS(\ov X)_{\rm tors}|$ and 
$|\NS(\ov Y)_{\rm tors}|$, we apply 
(\ref{e5}) to $X$, $Y$ and $X \times Y$, and obtain from Corollary \ref{kun}
and Section \ref{1.12} the decomposition of $\Ga$-modules
\begin{equation}\Br(\ov X\times\ov Y)_{\ell^m}=
\Br(\ov X)_{\ell^m}\oplus\Br(\ov Y)_{\ell^m} \oplus
\Hom(B^t_{\ell^m},A_{\ell^m})/\big(\Hom(\ov B^t,\ov A)/{\ell^m}\big).
\label{e2}\end{equation}
The case when $X$ and $Y$ are elliptic curves was considered in
\cite[Prop. 3.3]{SZCrelle}.
\end{sect}

\section{Proof of Theorem A}

The proof of Theorem A crucially uses the following
properties.
Let $C$ and $D$ be abelian varieties over a field $k$
finitely generated over its prime subfield, $\fchar(k)\not=2$.
Then

\smallskip

\noindent(1) the $\Gamma$-modules $V_{\ell}(C)$ and $V_{\ell}(D)$ are
semisimple, and the natural injective map
$$\Hom(C,D)\otimes \Q_{\ell} \hookrightarrow
\Hom_{\Gamma}(V_{\ell}(C),V_{\ell}(D)$$ is bijective;

\smallskip

\noindent(2) for almost all primes $\ell$
 the $\Gamma$-modules $C_{\ell}$ and
$D_{\ell}$ are semisimple, and the natural injective map
$$\Hom(C,D)/\ell \hookrightarrow
\Hom_{\Gamma}(C_{\ell},D_{\ell})$$ is bijective.
\medskip

Statement (1) was proved by the second named author in
characteristic $p>2$ \cite{ZarhinIz,ZarhinMZ1}, and by Faltings
\cite{Faltings1,Faltings2} in characteristic zero. Statement (2)
was proved by the second named author  in \cite[Thm.
1.1]{ZarhinMZ2} and \cite[Cor. 5.4.3 and Cor. 5.4.5]{ZarhinIn}, see also
\cite[Prop. 3.4]{SZJAG}, \cite[Thm. 4.4]{ZarhinMatSb2010},
\cite{ZarhinP} and \cite{ZarhinJIMJ}.

\smallskip

Theorem A is a consequence of the following result.

\begin{thm}
\label{main} Let $k$ be a field finitely generated over
its prime subfield. Let $X$ and $Y$ be smooth, projective
and geometrically integral varieties over
$k$. Then we have the following statements.

{\rm (i)} If  $\fchar(k)=0$, then
$[\Br(\ov X\times \ov Y)/(\Br(\ov X) \oplus
\Br(\ov Y))]^{\Gamma}$ is finite.

{\rm (ii)} If $\fchar(k)=p>2$, then the group
$[\Br(\ov X\times \ov Y)/(\Br(\ov X) \oplus
\Br(\ov Y))]^{\Gamma}({\non}p)$ is finite.
\end{thm}
\noindent{\bf Proof.}
Since $\Br(\ov X\times\ov Y)$ is a torsion group,
it is enough to prove these statements:

\smallskip

\noindent(a) {\it If $\ell$ is a prime, $\ell\not=\fchar(k)$, then
$V_\ell\big((\Br(\ov X\times \ov Y)/(\Br(\ov X) \oplus
\Br(\ov Y)))^{\Gamma}\big)=0$.}

\smallskip

\noindent(b) {\it For almost all primes $\ell$ we have
$\big(\Br(\ov X\times \ov Y)_\ell/(\Br(\ov X)_\ell\oplus\Br(\ov Y)_\ell)\big)^\Ga=0$.}

\medskip

Let us prove (a). Using (\ref{e1}) we obtain
\begin{eqnarray*}
V_\ell\big((\Br(\ov X\times\ov Y)/(\Br(\ov X)\oplus\Br(\ov Y)))^{\Gamma}\big)&=\\
V_\ell\big(\Br(\ov X\times\ov Y)/(\Br(\ov X) \oplus \Br(\ov Y))\big)^{\Gamma}&=\\
\big( V_\ell(\Br(\ov X\times\ov Y))/(V_\ell(\Br(\ov X))\oplus 
V_\ell(\Br(\ov Y)))\big)^{\Gamma}&=\\
\big(\Hom_{\Q_\ell}\big(V_{\ell}(B^t),V_{\ell}(A)\big)/\Hom(\ov{B^t},\ov{A})\otimes\Q_{\ell}\big)^{\Gamma}.&
\end{eqnarray*}
By a theorem of Chevalley
\cite[p. 88] {Chev} the semisimplicity of $\Gamma$-modules
$V_{\ell}(B^t)$ and $V_{\ell}(A)$ implies the semisimplicity of
the $\Gamma$-module
$\Hom_{\Q_{\ell}}(V_{\ell}(B^t),V_{\ell}(A))$.
From this we deduce
$$V_\ell\big((\Br(\ov X\times\ov Y)/(\Br(\ov X)\oplus\Br(\ov Y)))^{\Gamma}\big)=
\Hom_{\Gamma}(V_{\ell}(B^t),V_{\ell}(A))/\Hom(B^t,A)\otimes\Q_{\ell}=
0,$$ thus proving (a).

\medskip

Let us prove (b). By (\ref{e2}) it is enough to show that
$$\big(\Hom(B^t_\ell,A_\ell)/(\Hom(\ov B^t,\ov A)/\ell)\big)^{\Gamma}=0.$$
Since $\Hom(\ov {B^t},\ov{A})^{\Gamma}=\Hom(B^t,A)$,
the exact sequence
$$0\to\Hom(\ov{B^t},\ov{A})^{\Gamma}/\ell\to \big(\Hom(\ov{B^t},\ov{A})/\ell\big)^{\Gamma}
\to\H^1(k,\Hom(\ov{B^t},\ov{A}))$$
implies that for all but finitely many primes $\ell$ we have
$$\big(\Hom(\ov{B^t},\ov{A})/\ell\big)^{\Gamma}=\Hom(B^t,A)/\ell.$$
If we further assume that $\ell>2\dim(A)+2\dim(B)-2$, then,
by a theorem of Serre \cite{SerreTensor}, the
semisimplicity of the $\Gamma$-modules $B^t_{\ell}$ and
$A_{\ell}$ implies the semisimplicity of
$\Hom(B^t_{\ell},A_{\ell})$.
Hence we obtain
\begin{eqnarray*}
\big(\Hom(B^t_{\ell},A_{\ell})/(\Hom(\ov{B^t},\ov{A})/\ell)\big)^{\Gamma}&=&
\Hom(B^t_{\ell},A_{\ell})^{\Gamma}/(\Hom(\ov{B^t},\ov{A})/\ell)^{\Gamma}=\\
\Hom_{\Gamma}(B^t_{\ell},A_{\ell})/(\Hom(B^t,A)/\ell)&=&0,
\end{eqnarray*}
thus proving (b). QED

\begin{cor}
Let $k$ be a field finitely generated over
its prime subfield. Let $X$ and $Y$ be smooth, projective
and geometrically integral varieties over
$k$. Then we have the following statements.

{\rm (i)} Assume $\fchar(k)=0$. The group
$\Br(\ov X\times \ov Y)^\Gamma$ is finite if and only if
the groups $\Br(\ov X)^\Gamma$ and $\Br(\ov Y)^{\Gamma}$ are finite.

{\rm (ii)} Assume that $\fchar(k)$ is a prime $p>2$. The group
$\Br(\ov X\times \ov Y)^{\Gamma}({\non}p)$ is finite if and only if
the groups $\Br(\ov X)^{\Gamma}({\non}p)$ and
$\Br(\ov Y)^{\Gamma}({\non}p)$ are finite.
\end{cor}

\section{Proof of Theorem B} \label{3}

It is enough to prove the following statements:

\smallskip

\noindent (a) {\it The subgroup of $\Br(X\times Y)$ generated by
$\Br_1(X\times Y)$ and the images of $\Br(X)$ and $\Br(Y)$, has finite index.}

\noindent (b) {\it The cokernel of the natural map
$\Br_1(X)\oplus \Br_1(Y)\to \Br_1(X\times Y)$ is finite.}

\smallskip

Each of these statements formally follows from Theorem A, 
the functoriality of the Hochschild--Serre spectral sequence
\begin{equation}
E^{p,q}_2=\H^p(k,\H^q_\et(\ov X,\G_m))\Rightarrow \H^{p+q}_\et(X,\G_m) \label{ss}
\end{equation}
with respect to $X$, and the finiteness property stated in Remark \ref{1.9}.

Let us recall how (\ref{ss}) is usually applied.
If $X(k)\not=\emptyset$, then the canonical map  
$$E^{3,0}_2=\H^3(k,\bar k^*)\to \H^3_\et(X,\G_m)$$ 
has a section given by a $k$-point on $X$, and hence is injective.
The same is obviously true if $\H^3(k,\bar k^*)=0$.
The standard theory of spectral sequences now implies that 
the kernel of the canonical map
$$E^{0,2}_2=\Br(\ov X)^\Gamma\to E^{2,1}_2=\H^2(k,\Pic(\ov X))$$ is
the image of $\Br(X)$ in $\Br(\ov X)^\Gamma$.

\medskip

Let us prove (a). By functoriality of the spectral sequence
(\ref{ss}) we have the following commutative diagram with exact rows:
$$
\begin{array}{ccccc}
\Br(X\times Y)&\to&\Br(\ov X\times \ov Y)^\Gamma&\to &
\H^2(k,\Pic(\ov X\times\ov Y))\\
\uparrow&&\uparrow&&\uparrow\\
\Br(X)\oplus \Br(Y)&\to&\Br(\ov X)^\Gamma\oplus \Br(\ov Y)^\Gamma&\to &
\H^2(k,\Pic(\ov X))\oplus\H^2(k,\Pic(\ov Y))
\end{array}
$$
Note that the middle vertical map here is injective.
To prove (a) we must show that the image of
$\Br(X)\oplus\Br(Y)$ in $\Br(\ov X\times \ov Y)^\Gamma$
has finite index in the subgroup of the elements
that go to 0 in $\H^2(k,\Pic(\ov X\times\ov Y))$.
This follows from Theorem \ref{main} (i) and Remark \ref{1.9}.

To prove (b) we consider another commutative diagram with
exact rows, also constructed using the functoriality of
the spectral sequence (\ref{ss}):
$$
\begin{array}{ccccccc}
\Br(k)&\to&\Br_1(X\times Y)&\to &
\H^1(k,\Pic(\ov X\times\ov Y))&\to &0\\
&&\uparrow&&\uparrow&&\\
&&\Br_1(X)\oplus \Br_1(Y)&\to &
\H^1(k,\Pic(\ov X))\oplus\H^1(k,\Pic(\ov Y))&\to &0
\end{array}
$$
Statement (b) follows from this diagram and Remark \ref{1.9}.

\section{Proof of Theorem C}

The inclusion of the left hand side into the right hand side 
follows from functoriality of the Brauer group.
Thus we can assume that $X(\A_k)^\Br$ and $Y(\A_k)^\Br$ are not empty.
Since the Brauer group of a smooth projective variety is a torsion
group, to prove the opposite inclusion it is enough to show that
for any positive integer $n$ the subgroup
$\Br(X\times Y)_n$ is generated by the images of $\Br(X)_n$ and $\Br(Y)_n$,
together with some elements that pair trivially with 
$X(\A_k)^\Br \times Y(\A_k)^\Br$ with respect to the Brauer--Manin pairing.
The Kummer sequence gives a 
surjective map $\H^2_\et(X\times Y, \mu_n)\to \Br(X\times Y)_n$, so 
it suffices to show that, modulo the images of $\H^2_\et(X, \mu_n)$
and $\H^2_\et(Y, \mu_n)$, the group $\H^2_\et(X\times Y, \mu_n)$ is generated
by the elements that pair trivially with 
$X(\A_k)^\Br \times Y(\A_k)^\Br$. 

If $Z/X$ is a torsor under a $k$-group of multiplicative
type $G$ annihilated by $n$, then the type of $Z/X$, as recalled in 
Section \ref{ttype}, is 
the image of the class $[Z/X]$ under the composed map
$$\H^1_\et(X,G)\to \H^1_\et(\ov X,G)^\Ga\to 
\Hom_k(\hat G,\Pic(\ov X))=\Hom_k(\hat G,\Pic(\ov X)_n).$$
Recall that $S_X$ denotes the $k$-group scheme of multiplicative type 
that is dual to the $\Ga$-module $\Pic(\ov X)_n$. 

\ble
If $X(\A_k)^\Br$ is not empty, then there exists an $X$-torsor
under $S_X$ whose type is the identity map on $\hat S_X$.
\ele
\noindent{\bf Proof}. One of the main results of the descent theory 
of Colliot-Th\'el\`ene and Sansuc says that if 
$X(\A_k)^\Br\not=\emptyset$, then 
for any homomorphism of $\Ga$-modules $\hat G\to\Pic(\ov X)$
there exists an $X$-torsor under $G$
of this type, see \cite[Cor. 6.1.3 (1)]{S}. QED

\medskip

Let us choose one such $X$-torsor under $S_X$, and call it $\T_X$. 
(It is unique up to twisting by a $k$-torsor under $S_X$.)
Then $\ov \T_X$ is isomorphic to the $\ov X$-torsor $\T_{\ov X}$
from Section \ref{K}.
As in Remark \ref{rem14} we form the class 
$$[\T_X]\cup[\T_Y]\in\H^2_\et(X\times Y, S_X\otimes S_Y).$$
Pairing with it gives a map
$$\e:\Hom_k(S_X\otimes S_Y,\mu_n)=\Hom_k(S_X,\hat S_Y) \ \lra\ 
\H^2_\et(X\times Y,\mu_n).$$
For $\varphi\in\Hom_k(S_X,\hat S_Y)$ we can write
$\e(\varphi)=\varphi_*[\T_X]\cup [\T_Y]$, where
$\cup$ stands for the cup-product pairing
$$\H^1_\et(X,\hat S_Y)\times \H^1_\et(Y,S_Y)\to 
\H^2_\et(X\times Y,\mu_n).$$
Remark \ref{rem14} gives a commutative diagram
\begin{equation}
\begin{array}{ccc}
\Hom_k(S_X,\hat S_Y) &\stackrel{\e}\lra&\H^2_\et(X\times Y,\mu_n)\\
||&&\downarrow\\
\big(\H^1_\et(\ov X,\Z/n)\otimes\H^1_\et(\ov Y,\mu_n)\big)^\Ga
& \stackrel{\cup}\lra&\H^2_\et(\ov X\times \ov Y,\mu_n)^\Ga
\end{array}\label{com}\end{equation}
%By Theorem \ref{kun0} the homomorphism
%$$
%(\pi_X^*,\pi_Y^*,\e):\H^2_\et(\ov X,\mu_n)\oplus\H^2_\et(\ov Y,\mu_n)\oplus
%\Hom(S_X,\hat S_Y) \ \lra \ 
%\H^2_\et(\ov X\times \ov Y,\mu_n)
%$$
%is an isomorphism of $\Ga$-modules.

It is clear that Theorem C is a consequence of
Lemmas \ref{t2} and \ref{t3} below.

\ble \label{t2}
%Let $\pi_X:X\times Y\to X$ and $\pi_Y:X\times Y\to Y$ be natural projections.
We have 
$\H^2_\et(X\times Y,\mu_n)=\pi_X^*\H^2_\et(X,\mu_n)+\pi_Y^*\H^2_\et(Y,\mu_n)+\Im(\e)$.
\ele

\ble \label{t3}
For any positive integer $n$ we have the inclusion
$$X(\A_k)^{\Br_1(X)[n]}\times Y(\A_k)^{\Br_1(Y)[n]}\subset (X\times Y)(\A_k)^{\Im(\e)}.$$
\ele

\noindent{\bf Proof of Lemma \ref{t2}.} We use the spectral sequence
\begin{equation}
E^{p,q}_2=\H^p(k,\H^q_\et(\ov X,\mu_n))\Rightarrow \H^{p+q}_\et(X,\mu_n).
\label{ss2}
\end{equation}
Let us make a few observations in the case when $X$ is a smooth,
projective and geometrically integral variety over a number field $k$
such that $X(\A_k)\not=\emptyset$.
The canonical maps
$$E^{p,0}_2=\H^p(k,\mu_n)\to \H^p_\et(X,\mu_n)$$
are injective for $p\geq 3$. Indeed, for such $p$ the natural map
$$\H^p(k,M)\to \bigoplus_{k_v\simeq\R} \H^p(k_v,M)$$
is a bijection for any finite Galois module $M$, see \cite[Thm. I.4.10 (c)]{ADT}.
Next, the natural map $\H^p(k_v,M)\to \H^p_\et(X\times_k k_v,M)$
is injective for any $p$ 
since any $k_v$-point of $X$ defines a section of it.
It follows that the composite map 
$$\H^p(k,M)\to\H^p_\et(X,M)\to\bigoplus_{k_v\simeq\R}\H^p_\et(X\times_k k_v,M)$$
is injective, and this implies our claim. 
We note that 
$$E^{2,0}_2=\H^2(k,\mu_n)\to \H^2_\et(X,\mu_n)$$
is also injective. The argument is similar once we identify
$\H^2(k,\mu_n)=\Br(k)_n$ using the Kummer sequence and Hilbert's
theorem 90, and use the
embedding of $\Br(k)$ into the direct sum of $\Br(k_v)$, for
all completions $k_v$ of $k$, provided by global class field theory,
together with the existence of $k_v$-points on $X$ for every place $v$.
This implies the triviality of all the canonical maps in the spectral sequence
whose target is $E^{p,0}_2=\H^p(k,\mu_n)$ for $p\geq 2$.

Let us write $\tilde\H^2_\et(X,\mu_n)$ 
for the quotient of $\H^2_\et(X,\mu_n)$ by the (injective) image of $\H^2(k,\mu_n)$.
Using the above remarks we obtain from (\ref{ss2}) the following exact
sequence:
\begin{equation}
0\to \H^1(k,\H^1_\et(\ov X,\mu_n))\to 
\tilde\H^2_\et(X,\mu_n) \to 
\H^2_\et(\ov X,\mu_n)^\Ga \to
\H^2(k,\H^1_\et(\ov X,\mu_n)).\label{ee4}
\end{equation}
There are similar sequences for $Y$ and $X\times Y$ linked by the
maps $\pi_X^*$ and $\pi_Y^*$.

Let us define 
$$\sH=\pi_X^*\H^2_\et(X,\mu_n)+\pi_Y^*\H^2_\et(Y,\mu_n)+\Im(\e)\subset 
\H^2_\et(X\times Y,\mu_n).$$
It is clear that the (injective) image of $\H^2(k,\mu_n)$ in
$\H^2_\et(X\times Y,\mu_n)$ is contained in $\sH$, 
so to prove Theorem C it is enough to prove that the natural map
$\sH\to\tilde\H^2_\et(X\times Y,\mu_n)$ is surjective.

By (\ref{Ku1}) the image of 
$\H^1(k,\H^1_\et(\ov X\times \ov Y,\mu_n))\to \tilde\H^2_\et(X\times Y,\mu_n)$
is contained in $\sH$. In view of (\ref{ee4}) it remains to show that
every element of the kernel of the map
$$\H^2_\et(\ov X\times \ov Y,\mu_n)^\Ga \to 
\H^2(k,\H^1_\et(\ov X\times \ov Y,\mu_n))$$
comes from $\sH$.
By Theorem \ref{kun0} and (\ref{Ku1}) this map can be written as
$$\H^2_\et(\ov X,\mu_n)^\Ga\oplus \H^2_\et(\ov Y,\mu_n)^\Ga \oplus
\Hom_k(S_X,\hat S_Y)\lra\H^2(k,\H^1_\et(\ov X,\mu_n))\oplus 
\H^2(k,\H^1_\et(\ov Y,\mu_n)).$$
By the commutativity of the diagram (\ref{com}) for any 
$\varphi\in \Hom_k(S_X,\hat S_Y)$, the element
$\e(\varphi)\in \H^2_\et(X\times Y,\mu_n)$ maps to 
$$\varphi\in \Hom(S_X,\hat S_Y)^\Ga\subset\H^2_\et(\ov X\times \ov Y,\mu_n)^\Ga.$$
This implies that
for any $a\in\H^2_\et(X\times Y,\mu_n)$ there 
exists an element $b\in\sH$ such that the image of $a-b$ in 
$\H^2_\et(\ov X\times \ov Y,\mu_n)^\Ga$ is $\pi_X^*(x)+\pi_Y^*(y)$ for some
$x\in \H^2_\et(\ov X,\mu_n)^\Ga$ and $y\in \H^2_\et(\ov Y,\mu_n)^\Ga$.
From the exact sequence (\ref{ee4}) for $X\times Y$ we see that
$\pi_X^*(x)+\pi_Y^*(y)$ goes to zero in 
$\H^2(k,\H^1_\et(\ov X,\mu_n))\oplus 
\H^2(k,\H^1_\et(\ov Y,\mu_n))$, hence $x$ goes to zero in 
$\H^2(k,\H^1_\et(\ov X,\mu_n))$, and $y$ goes to zero in 
$\H^2(k,\H^1_\et(\ov Y,\mu_n))$. By (\ref{ee4}) for $X$
we see that $x$ is the image of some $c\in \H^2_\et(X,\mu_n)$.
Similarly, $y$ is the image of some $d\in \H^2_\et(Y,\mu_n)$.
This proves that $a-(b+\pi_X^*(c)+\pi_Y^*(d))$ goes to zero in 
$\H^2_\et(\ov X\times \ov Y,\mu_n)^\Ga$, and hence belongs to $\sH$.
Thus $a\in\sH$. QED

\medskip

\noindent{\bf Proof of Lemma \ref{t3}.} 
Let $M$ be a finite $\Ga$-module such that $nM=0$.
Let $M^D$ be the dual module $\Hom(M,\bar k^*)$. 
If $v$ is a non-archimedean place of $k$, we write
$\H^1_{\rm nr}(k_v,M)$ for the unramified subgroup of 
$\H^1(k_v,M)$. By definition, it consists of the classes
that are annihilated by the restriction to the maximal
unramified extension of $k_v$. We write $P^1(k,M)$
for the restricted product of the abelian groups
$\H^1(k_v,M)$ relative to the subgroups
$\H^1_{\rm nr}(k_v,M)$, where $v$ is a non-archimedean place of $k$.
By \cite[Lemma I.4.8]{ADT} the image of the diagonal map
$$\H^1(k,M)\to \prod_{{\rm all}\ v}\H^1(k_v, M)$$
is contained in $P^1(k,M)$. Let us denote this image by $U^1(k,M)$.

The local pairings 
$\H^1(k_v,M)\times \H^1(k_v,M^D)\to \H^2(k_v,\mu_n)$
give rise to the
global Poitou--Tate pairing 
$$(\ ,\ ):\quad P^1(k,M)\times P^1(k,M^D)\ \lra \ \Z/n.$$
It is a perfect duality of locally compact abelian groups,
moreover, the subgroups $U^1(k,M)$ and $U^1(k,M^D)$
are exact annihilators of each other \cite[Thm. I.4.10 (b)]{ADT}.

Let $\varphi\in\Hom_k(S_X,\hat S_Y)$. 
Let $(P_v)\in X(\A_k)$ be an adelic point that is Brauer--Manin
orthogonal to $\Br_1(X)[n]$, and let $(Q_v)\in Y(\A_k)$ be an
adelic point orthogonal to $\Br_1(Y)[n]$. The Brauer--Manin
pairing of the adelic point $(P_v\times Q_v)$ with the image
of $\e(\varphi)=\varphi_*[\T_X]\cup [\T_Y]$ in $\Br(X\times Y)$
is given by the Poitou--Tate pairing, so
to prove Lemma \ref{t3} we need to show that 
\begin{equation}
(\varphi_*[\T_X](P_v),[\T_Y](Q_v))=0, \label{ee1}
\end{equation}
where in the above notation $M=\hat S_Y$, $M^D=S_Y$.
We point out that for any $a\in \H^1(k,\hat S_Y)$ we have
$a\cup [\T_Y]\in \Br_1(Y)[n]$, and hence $(a,[\T_Y](Q_v))=0$.
If an element of $P^1(k,S_Y)$ is orthogonal to $U^1(k,\hat S_Y)$, 
then it belongs to $U^1(k,S_Y)$. Therefore, we must have 
\begin{equation}
[\T_Y](Q_v)\in U^1(k,S_Y). \label{ee2}
\end{equation}
Similarly, for any $b\in \H^1(k,S_Y)$ we have
$\varphi_*[\T_X]\cup b\in\Br_1(X)[n]$, and hence
$(\varphi_*[\T_X](P_v),b)=0$. Since every element of $P^1(k,\hat S_Y)$ 
orthogonal to $U^1(k,S_Y)$ belongs to $U^1(k,\hat S_Y)$,
this implies 
\begin{equation}
\varphi_*[\T_X](P_v)\in U^1(k,\hat S_Y). \label{ee3}
\end{equation}
Now (\ref{ee2}) and (\ref{ee3}) imply (\ref{ee1}). QED

\end{document}